\documentclass[12pt]{article}
\usepackage{graphicx}
\usepackage{amssymb}
\usepackage{amsmath}
\numberwithin{equation}{section}

\begin{document}
\author{Lev Sakhnovich}
\date{April 3, 2008}
\textbf{On the solutions of Knizhnik-Zamolodchikov system}

\begin{center} Lev Sakhnovich  \end{center}
735 Crawford ave., Brooklyn, 11223, New York, USA.\\
 E-mail address:lev.sakhnovich@verizon.net
\begin{center}Abstract \end{center}
We consider the Knizhnik-Zamolodchikov system of linear differential
equations. The coefficients of this system are rational functions.
We prove that under some conditions the solution of KZ system is
rational too.  We give the method of constructing the corresponding
rational solution. We deduce the asymptotic formulas for the
solution of KZ system when $\rho$ is integer.
\newpage
\section{Introduction }
1.We consider the differential system
\begin{equation}
\frac{{\partial}W(z_{1},z_{2},...,z_{n})}{{\partial}z_{j}}={\rho}A_{j}(z_{1},z_{2},...z_{n})W,
\quad  1{\leq}j{\leq}n,\end{equation} where
 $A_{j}(z_{1},z_{2},...,z_{n})$ and $W(z_{1},z_{2},...,z_{n})$ are $n{\times}n$ matrix
functions. We suppose that $A_{j}(z_{1},z_{2},...,z_{n})$ has the
form
\begin{equation}
A_{j}(z_{1},z_{2},...,z_{n})=\sum_{k=1,
k{\ne}j}^{n}\frac{P_{j,k}}{z_{j}-z_{k}},\end{equation} where
$z_{k}{\ne}z_{\ell}$ if $k{\ne}\ell$. Here the matrices $P_{j,k}$
are connected with the matrix representation of the symmetric group
$S_{n}$ and are defined by formulas (2.1)-(2.4). We note that the
well-known Knizhnik- Zamolodchikov equation has the form (1.1),(1.2)
(see [3]). This system has found applications in several  areas of
mathematics and physics (see [3],[5]). In the first part of the
paper (section 2) we prove the following assertion:

\emph{ The fundamental solution of KZ system (1.1),(1.2) is
rational,when ${\rho}$ is integer.}

We give an effective method of constructing this rational solution.
Our method is elementary and is based on the results of linear
algebra. The more complicated approach to constructing rational
solutions is given by  G.Felder and A.Veselov [4]. In our previous
papers [7],[9],[10] we constructed rational solutions only for the
cases $\rho=\pm{1}$. In sections 3-6 we use the variables which were
introduced by A.Varchenko [13]:
\begin{equation}u_{1}=z_{1}-z_{2},
u_{2}=\frac{z_{2}-z_{3}}{z_{1}-z_{2}},...,u_{n-1}=\frac{z_{n-1}-z_{n}}{z_{n-2}-z_{n-1}},
u_{n}=z_{1}+z_{2}+...+z_{n}.
\end{equation} In terms of the variables $u_{j}$ the system KZ takes
the form
\begin{equation}
\frac{{\partial}W(u_{1},u_{2},...,u_{n})}{{\partial}u_{j}}={\rho}H_{j}(u_{1},u_{2},...u_{n})W,
\quad  1{\leq}j{\leq}n.\end{equation}We investigate the form of
$H_{j}(u_{1},u_{2},...u_{n})$ in a more detailed manner than it is
done in paper [13]. In particular we prove that
$H_{1}(u_{1},u_{2},...u_{n})$ depends only on $u_{1}$ and
$H_{1}(u_{1},u_{2},...u_{n}), (2{\leq}j{\leq}n)$ does not depend on
$u_{1}$. Using the results of sections 1-4 we deduce the asymptotic
formula for the solutions of KZ system (1.1) when $\rho$ is integer
and $u_{j}{\to}0.$ The corresponding result for the irrational
$\rho$ is well-known [13]. In the last section of the paper we
consider the examples (n=3 and n=4). As a by-product we obtain the
following result:\\ \emph{The hypergeometric function F(a,b,c) is
rational when}
\begin{equation}a=-\rho,\quad b=-3\rho,\quad c=1-2\rho,\quad (\rho\quad is \quad integer).\end{equation}
\section{Rational solutions of KZ system}
1. We consider the natural representation of the symmetric group
$S_{n}$ (see [2]).  By $(i;j)$ we denote the permutation which
transposes $i$ and $j$ and preserves all the rest. The $ n{\times}n$
matrix which corresponds to $(i;j)$ is denoted by
\begin{equation}
P(i,j)=\{p_{k,\ell}(i,j)\}_{k,\ell=1}^{n}.\quad (i{\ne}j).
\end{equation}The elements $\{p_{k,\ell}(i,j)\}$ are equal to zero
except
\begin{equation}\{p_{k,\ell}(i,j)\}=1, \quad when\quad either \quad
k=i,\ell=j \quad or \quad k=j, \ell=i ,\end{equation}
\begin{equation}\{p_{k,k}(i,j)\}=1, \quad when \quad k{\ne}i,j.\end{equation}
2. We begin with the first equation of the KZ system (1.1), (1.2):
\begin{equation}\frac{{\partial}W(z)}{{\partial}z_{1}}={\rho}A_{1}(z)W(z),\end{equation}
where $z=[z_{1},z_{2},...,z_{n}].$ In this section we assume that
$\rho$ is integer.\\
\textbf{Theorem 2.1.}\emph{Let $\rho$ be integer. Then the
fundamental matrix solution $W(z)$ of system (2.4) can be written in
the form}\begin{equation}
W(z)=\sum_{k=2}^{n}\sum_{j=1}^{m}\frac{L_{k,j}(\xi)}{(z_{1}-z_{k})^{j}}+Q(z),\end{equation}
\emph{where $\xi=[z_{2},...,z_{n}],\quad \mathrm{m}=|\rho|$ and
$L_{k,j}(\xi)$ are  $n{\times}n$  matrix functions, Q(z) is a matrix
polynomial in respect to $z_{1}$ and}
\begin{equation}\mathrm{deg}Q(z)=m(n-1),\quad if\quad \rho>0;\quad
\mathrm{deg}Q(z)=m,\quad if \quad\rho<0.\end{equation} \emph{Proof.}
Changing the variables $z_{1}=1/x$ we obtain
\begin{equation}\frac{{\partial}V}{{\partial}u}=B(u,\xi)V,\end{equation}
where \begin{equation}V(u,\xi)=W(1/u,\xi),\quad
B(u,\xi)=(-{\rho}/u)\sum_{k=2}^{n}P_{1,k}/(1-uz_{k}).\end{equation}
From relation (2.8) we deduce the representation
\begin{equation}
B(u,\xi)=-\rho\sum_{p=-1}^{\infty}u^{p}T_{p},\end{equation} where
\begin{equation}T_{p}=\sum_{k=2}^{n}P_{1,k}z_{k}^{p+1},\quad
p{\geq}-1.\end{equation}We investigate the case when $V(u.\xi)$ can
be written in the form
\begin{equation}V(u,\xi)=\sum_{j=s}^{\infty}u^{j}G_{j}(\xi),\quad
G_{s}{\ne}0,\quad |u|<r_{0}.\end{equation} Here $G_{j}(\xi)$ are
$n{\times}n$ matrix functions. It follows from (2.7), (2.9) and
(2.11) that
\begin{equation}[(q+1)+{\rho}T_{-1}]G_{q+1}={-\rho}\sum_{j+\ell=q}T_{j}G_{\ell},\quad
j{\geq}0.\end{equation} According to (2.10) the matrix $T_{-1}$ has
the following structure
\begin{equation}T_{-1}=(n-2)I_{n}+S , \end{equation}
where $n{\times}n$ matrix $S$ is defined by the equality
\begin{equation}S=\left[\begin{array}{cc}
                          2-n & e \\
                          e^{\tau}    & 0
                        \end{array}\right],\quad
                        e=\overbrace{[1,1,...,1]}^{n-1}.\end{equation}
Here $e^{\tau}$ is transposed $e$,i.e.
$e^{\tau}=\mathrm{col}[1,1,...,1]$. Let us write the eigenvalues of
$T_{-1}$:
\begin{equation}\lambda_{-1}=n-1,\quad \lambda_{-2}=n-2,\quad
\lambda_{-3}=-1.\end{equation}The corresponding $1{\times}n$
eigenvectors of $T_{-1}$ have the forms \begin{equation}
U_{1}=\mathrm{col}\overbrace{[1,1,...,1]}^n,\quad
U_{2,j}=\mathrm{col}[0,a_{1,j},a_{2,j},...,a_{n-1,j}],\end{equation}
\begin{equation}U_{3}=\mathrm{col}\overbrace{[n-1,-1,-1,...-1]}^n.\end{equation}
We note, that vectors $U_{2,j}, (1{\leq}j{\leq}n-2)$ are linearly
independent and satisfy the condition
\begin{equation}a_{1,j}+a_{2,j}+...+a_{n-1,j}=0.\end{equation}
 We shall construct $n$ linear independent
solutions $Y_{k}(z), (1{\leq}k{\leq}n)$ of system (2.4). We begin
with
the case $m=-\rho>0.$\\
\emph{Step 1.}We assume that
\begin{equation}G_{k}=0,\quad k{\leq}m(n-2)-1.\end{equation}
In this case we have (see(2.12))
\begin{equation}[(n-2)I_{n}-T_{-1}]G_{m(n-2)}=0.\end{equation}
According to (2.16)  the vector $G_{m(n-2)}(j)=U_{2,j}$ satisfies
the relation (2.20). From formula (2.12) we find the coefficients
$G_{p}(j),( m(n-2)<p<m(n-1)).$ In case $p=m(n-1)$ the matrix
$[(n-1)I_{n}-T_{-1}]$ is not invertible because the equality
\begin{equation} [(n-1)I_{n}-T_{-1}]U_{1}=0 \end{equation}
is true. Now we shall use the relations
\begin{equation}P_{1,k}U_{1}=U_{1},\quad 2{\leq}k{\leq}n;
\quad (U_{2,j},U_{1})=0,\quad 1{\leq}j{\leq}n-2.\end{equation} It
follows from (2.22) that the right side of equality (2.12),when
$q+1=m(n-1),$ is orthogonal to $U_{1}$. Hence the corresponding
equation (2.12), when  $q+1=m(n-1),$ has a solution $G_{m(n-1)}(j)$.
To find the matrix coefficients $L_{k,p}(j)$ we consider the block
matrices
\begin{equation}L_{k}(j)=\mathrm{col}[L_{k,1}(j), L_{k,2}(j), ...,
L_{k,m}(j)],\quad
2{\leq}k{\leq}n,\end{equation}\begin{equation}L(j)=\mathrm{col}[L_{2}(j),
L_{3}(j), ..., L_{n}(j)].\end{equation}
\begin{equation}G(j)=\mathrm{col}[G_{1}(j),
G_{2}(j), ..., G_{m(n-1)}(j)],\end{equation}and use  the equality
\begin{equation} RL_{j}=G_{j} .\end{equation}  The
$mn(n-1){\times}mn(n-1)$ matrix $R$ has the following form:
\begin{equation}R=[R_{1}, R_{2},...,R_{n-1}], \end{equation}
where the $n{\times}n$ blocks $R_{k}$ of the $mn(n-1){\times}mn$
matrix $R$ are defined by the relations
\begin{equation} \{R_{k}\}_{p,s}=0_{n,n},\quad p<s;\quad
\{R_{k}\}_{p,s}=z_{k}^{s}\left(\begin{array}{c}
                           p-1 \\
                           s-1
                         \end{array}\right)
I_{n}\quad p{\geq}s.\end{equation} Here $0_{n,n}$ is the
$n{\times}n$ zero matrix, $1{\leq}k{\leq}n-1, 1{\leq}p{\leq}m(n-1),
1{\leq}s{\leq}m.$ The introduced matrix $R$ is the block confluent
Vandermonde matrix. Hence the matrix $R$ is invertible and we have
\begin{equation}L_{j}=R^{-1}G_{j}.\end{equation}In paper [] we
proved the assertion(necessary condition)\\\emph{If the solution
$Y_{j}(z)$ of system (2.4) is rational then relation
(2.26) is true.}\\
 In such a way we constructed $(n-2)$ linearly
independent vector functions $Y_{j}(z) \quad (1{\leq}j{\leq}n-2)$ of
the form \begin{equation}
Y_{j}(z)=\sum_{k=2}^{n}\sum_{p=1}^{m}\frac{L_{k,p}(j)(\xi)}{(z_{1}-z_{k})^{p}}.\end{equation}
 Later we
shall show that $Y_{j}$ are solutions of system (2.4).\\
\emph{Step 2.} Now we assume that
\begin{equation}G_{k}=0,\quad -m{\leq}k{\leq}m(n-1)-1.\end{equation}
In this case we have
\begin{equation}G_{m(n-1)}=U_{1}.\end{equation} It follows from the
relation
\begin{equation}P_{1,k}U_{1}=U_{1} \end{equation}
that the vector function
\begin{equation}Y_{n-1}(z)=\prod_{k=2}^{n}(z_{1}-z_{k})^{-m}U_{1}
\end{equation} is a solution of system (2.4)\\
\emph{Step 3.} To construct the solution $Y_{n}(z)$ we consider the
case when
\begin{equation}G_{-m}=U_{3}.\end{equation}
Using relations (2.12) we find the coefficients $G_{p},
(-m{\leq}p{\leq}m(n-2)-1).$ In the case $p=m(n-2)$ the matrix
$[m(n-2)I_{n}-mT_{-1}]$ is not invertible. We represent the right
side of (2.12) when $q+1=m(n-2)$ in the form $U_{2}+{\beta}U_{3}$,
where  $U_{2}$ is the  fixed concrete vector .We remark that the
vector $U_{2}$ is a linear combination of the vectors $U_{2}(j)$.
The vector $U_{3}$ is orthogonal to $U_{2}(j)$. Hence the equation
\begin{equation} [m(n-2)I_{n}-mT_{-1}]G_{m(n-2)}={\beta}U_{3}
\end{equation} has a solution $G_{m(n-2)}$. We consider the case
when \begin{equation}L_{n,j}=0,\quad
(1{\leq}j{\leq}m).\end{equation} We note that in case (2.35)
equality (2.29) takes the form
\begin{equation}L=R^{-1}G,\end{equation}where
\begin{equation}
G=\mathrm{col}[G_{1},G_{2},...,G_{m(n-2)}],\quad
R=[R_{1},R_{2},...,R_{n-2}].\end{equation} By relations (2.12) and
(2.38) we find $G_{p}\quad(-m{\leq}p{\leq}0)$ and $L_{k,j},\quad
(2{\leq}k{\leq}n-1,\quad 1{\leq}j{\leq}m).$ In such a way we
construct the vector function
\begin{equation}Y_{n}(z)=\sum_{k=2}^{n-1}\sum_{p=1}^{m}\frac{L_{k,p}(\xi)}{(z_{1}-z_{k})^{p}}+Q(z_{1},\xi).\end{equation}
where $Q(z_{1},\xi)$ is a polynomial  of degree $m$ in respect to
$z_{1}$.\\ \emph{Step 4.} To prove that the constructed vector
functions $Y_{k}(z), (1{\leq}k{\leq}n-2)$ are solutions of system
(2.4) we consider the vector functions
\begin{equation}
M_{k}(z_{1},\xi)=\frac{{\partial}Y_{k}}{{\partial}z_{1}}-A_{1}(z_{1},\xi)Y_{k}.\end{equation}
The vector functions $M_{k}(z_{1},\xi)$ are rational in respect to
$z_{1}$ with the poles in the points $z_{1}=z_{j},
(2{\leq}j{\leq}n)$. It follows from (2.29) that
\begin{equation}M_{k}(z_{1},\xi)=O(|z_{1}|^{-p}),\quad
z_{1}{\to}\infty,\end{equation} where $p=m(n-1)+2.$ We introduce the
polynomials
\begin{equation}q(z_{1},\xi)=\prod_{s=2}^{n}(z_{1}-z_{s})^{m},\quad
q_{j}(z_{1},\xi)=\prod_{s=2,
 s{\ne}j}^{n}(z_{1}-z_{s})^{m}.\end{equation} Using relations (2.42)
and (2.43)  we deduce that
\begin{equation}q(z_{1},\xi)M_{k}(z_{1},\xi)=R_{k}(z_{1},\xi)+\sum_{k=2}^{n}\frac{(-m)(I-P_{1,k})}{z_{1}-z_{k}}
L_{k,m}c_{k}
\end{equation}
Here $R_{k}(z_{1},\xi)$ is a polynomial with respect to $z_{1}$ and
\begin{equation}
c_{k}=q_{k}(z_{k},\xi)=\prod_{2{\leq}s{\leq}n,s{\ne}k}(z_{k}-z_{s})^{m}.\end{equation}
 Comparing (2.42) and (2.44) we have
\begin{equation}R_{k}(z_{1},\xi)=0,\quad
\sum_{k=2}^{n}(I-P_{1,k}) L_{k,m}c_{k}=0.
\end{equation}By $X_{k}$ we denote the vectors such that
\begin{equation} P_{1,k}X_{k}=-X_{k},\quad X_{k}{\ne}0,\quad
2{\leq}k{\leq}n.\end{equation} From (2.47) we infer that
\begin{equation}(I_{n}-P_{1,k})L_{k,m}={\alpha}_{k}X_{k}.\end{equation}
The vectors $X_{k} \quad (2{\leq}k{\leq}n)$ are linearly
independent. Hence according to (2.46)  we have
\begin{equation}(I_{n}-P_{1,k})L_{k,m}=0.\end{equation}
It follows directly from (2.44), (2.46) and (2.49) that
$M_{k}(z_{1},\xi)=0 \quad (2{\leq}k{\leq}n).$ Consequently, the
vector functions $Y_{k}(z)
\quad (2{\leq}k{\leq}n)$ are solutions of system (2.4).\\
It was shown before that the vector function $Y_{n-1}(z)$  is a
solution of system (2.4).\\
\emph{Step 5.} Now we shall prove that the vector function
$Y_{n}(z)$ (see (2.40))is a solution of system (2.4).To do it we
introduce the rational vector function
\begin{equation}
M_{n}(z_{1},\xi)=\frac{{\partial}Y_{n}}{{\partial}z_{1}}-A_{1}(z_{1},\xi)Y_{n}.\end{equation}
From (2.37) and (2.38) we infer the equality
\begin{equation}z_{1}^{m(n-2)+1}M_{n}(z_{1},\xi){\to}\tilde{u},\quad
when \quad z_{1}{\to}\infty.\end{equation}Here vector $\tilde{u}$ is
a fixed vector which satisfies condition (2.18). We introduce the
polynomials
\begin{equation}r(z_{1},\xi)=\prod_{s=2}^{n-1}(z_{1}-z_{s})^{m},\quad
r_{j}(z_{1},\xi)=\prod_{s=2,
 s{\ne}j}^{n-1}(z_{1}-z_{s})^{m}.\end{equation}
In view of (2.50) and (2.52) the equality
\begin{equation}r(z_{1},\xi)M_{n}(z_{1},\xi)=R_{n}(z_{1},\xi)+\sum_{k=2}^{n-1}\frac{(-m)(I-P_{1,k})}{z_{1}-z_{k}}
L_{k,m}d_{k}
\end{equation} is true. Here $R_{k}(z_{1},\xi)$ is a polynomial with respect to $z_{1}$ and
\begin{equation}
d_{k}=r_{k}(z_{k},\xi)=\prod_{2{\leq}s{\leq}n-1,s{\ne}k}(z_{k}-z_{s})^{m}.\end{equation}
According to (2.51) the relations
\begin{equation}R_{n}(z_{1},\xi)=0,\quad
\sum_{k=2}^{n-1}(-m)(I-P_{1,k})L_{k,m}d_{k} =\tilde{u}
\end{equation}hold. From (2.48) and (2.55) we infer
\begin{equation}\sum_{k=2}^{n-1}\alpha_{k}d_{k}X_{k}=\tilde{u}.\end{equation}
We took into account that $X_{1}$ is orthogonal to $\tilde{u}$. If
we consider the case $L_{p,j}=0,\quad (1{\leq}j{\leq}m,\quad
p{\ne}n)$, we obtain the analogue of (2.56):
\begin{equation}\sum_{k=2, k{\ne}p}^{n}\tilde{\alpha_{k}}\tilde{d_{k}}X_{k}=\tilde{u}.\end{equation}
Using the independence of the system $X_{k}\quad (2{\leq}k{\leq}n)$
and comparing (2.56) and (2.57) we receive the equalities
${\alpha}_{p}=0,\quad (2{\leq}p{\leq}n).$ Hence we have
\begin{equation}\tilde{u}=0.\end{equation}.It follows from (2.53), (2.55) and (2.58) that
$M_{n}(z_{1},\xi)=0$ and the vector function $Y_{n}(z)$ is a
solution of system (2.4). It is easy to see that  the constructed
solutions $Y_{k}(z)$ are linearly independent.\\Now we use the
following proposition [10]:\\ \emph{If $W(z)$ is a rational
fundamental solution of system (2.4),when $\rho=-m$, then
$[W^{-1}(z)]^{\tau}$ is a fundamental rational solution of (2.4),
when $\rho=m$.}\\Hence the assertion of the theorem is true for all
integer $\rho.$ We note that relations (2.6) are true (see[8])
 The
theorem is proved.\\
Further we need the following assertion.\\
\textbf{Proposition 2.1.} \emph{If the matrix function}
$W_{j}(z_{1},...,z_{j},z_{j+1},...,z_{n})$\emph{ is a solution of
the equation}
\begin{equation}\frac{\partial}{{\partial}z_{j}}W_{j}(z)={\rho}A_{j}(z)W_{j}(z),\quad
(1{\leq}j{\leq}n)\end{equation}\emph{then the matrix function }
$W_{j+1}(z)=P_{j,j+1}W(z_{1},...,z_{j+1},z_{j},...,z_{n})$ \emph{is
a solution of the equation}
\begin{equation}\frac{\partial}{{\partial}z_{j+1}}W_{j+1}(z)={\rho}A_{j+1}(z)W_{j+1}(z),\quad
(1{\leq}j{\leq}n-1)\end{equation} \emph{Proof.} The proposition
follows directly from (2.1-(2.3)) and the relations
\begin{equation} P_{j,j+1}P_{j,i}P_{j,j+1}=P_{j+1,i};\quad
j{\ne}i,\quad j{\ne}i+1,\end{equation}
\begin{equation}P_{j,j+1}P_{j,j+1}P_{j,j+1}=P_{j+1,j}.\end{equation}
\textbf{Corollary 2.1.} \emph{Every equation of system (1.1) has a
fundamental rational solution when $\rho$ is integer.}\\
We introduce the fundamental $n{\times}n$ matrix solution of system
(2.4):
\begin{equation}W(z)=[Y_{1}(z),Y_{2}(z),...,Y_{n}(z)].\end{equation}
Now we fix some point $z_{1,0}$ and consider the new fundamental
matrix solution
\begin{equation}W_{1}(z)=W(z)W^{-1}(z_{1,0},\xi).\end{equation}
of system (2.4), which satisfies the condition
\begin{equation}W_{1}(z_{1, 0},\xi)=I_{n}.\end{equation}
In the next section we shall use the normalized fundamental solution
$W_{1}(z)$.
\section{Common solution, consistency}
1.In section 2 we have considered only one equation (2,4) of  KZ
system (1.1). Using this result we construct a common rational
solution of KZ system (1.1). We note that the KZ system is
consistent , i.e. the relations
\begin{equation}\rho(\frac{{\partial}A_{i}}{{\partial}z_{j}}-\frac{{\partial}A_{j}}{{\partial}z_{i}})+
{\rho}^{2}[A_{i},A_{j}]=0 \end{equation} are valid.\\
The correctness of (3.1) follows from the properties of $P_{i,j}$:\\
\begin{equation} P_{i,j}=P_{j,i},\end{equation}
\begin{equation} [P_{i,j}+P_{j,k},P_{i,k}]=0;\quad (i,j,k \quad are \quad distinct),\end{equation}
\begin{equation} [P_{i,j},P_{k,\ell}]=0;\quad (i,j,k,{\ell} \quad
are \quad distinct),\end{equation} The following assertion can be
easily verified (see [6],Ch.12.).\\
\textbf{Proposition 3.1.} \emph{Let the $n{\times}n$ matrix function
$W_{1}(z_{1},z_{2})$ satisfy the equation}
\begin{equation}\frac{\partial}{{\partial}z_{1}}W_{1}=B_{1}(z_{1},z_{2})W_{1}
\end{equation}\emph{and the condition}
\begin{equation}W_{1}(z_{1,0},z_{2})=I_{n},\end{equation}\emph{where
$z_{1,0}$ is a fixed point. If the relation}
\begin{equation}\frac{{\partial}B_{1}}{{\partial}z_{2}}-\frac{{\partial}B_{2}}{{\partial}z_{1}}+
[B_{1},B_{2}]=0 \end{equation}\emph{ is valid, then the matrix
function }
\begin{equation}U(z_{1},z_{2})=\frac{{\partial}W_{1}}{{\partial}z_{2}} - B_{2}(z_{1},z_{2})W_{1}
\end{equation}\emph{satisfies equation.} (3.5)\\
From  condition (3.6) and equality (3.8) we deduce the relation
\begin{equation}U(z_{1,0},z_{2})= -
B_{2}(z_{1,0},z_{2}).\end{equation} Using Proposition 3.1 and
relation (3.9) we have -
\begin{equation}\frac{{\partial}W_{1}}{{\partial}z_{2}}-
B_{2}(z_{1},z_{2})W_{1}=-W_{1}(z)B_{2}(z_{1,0},z_{2}).\end{equation}
Now we introduce the $n{\times}n$ matrix function
$W_{2}(z_{1,0},z_{2})$ which satisfies the equation
\begin{equation}\frac{{\partial}W_{2}}{{\partial}z_{2}}-B_{2}(z_{1,0},z_{2})W_{2}=0
\end{equation}
and  the condition
\begin{equation}W_{2}(z_{1, 0},z_{2,0})=I_{n},\end{equation}
where $z_{1,0}$ and $z_{2,0}$ are fixed points. It follows from
(3.11)
and (3.12) the assertion.\\
\textbf{Proposition 3.2.} \emph{ The $n{\times}n$ matrix funtion}
\begin{equation}W(z_{1},z_{2})=W_{1}(z_{1},z_{2})W_{2}(z_{1,0},z_{2})\end{equation}
\emph{satisfies the equations}
\begin{equation}\frac{{\partial}W}{{\partial}z_{k}}-B_{k}(z_{1},z_{2})W=0;\quad
k=1,2
\end{equation}
\emph{and  the condition}
\begin{equation}W(z_{1, 0},z_{2,0})=I_{n},\end{equation} We introduce the notation
\begin{equation}\mu_{i}=(z_{1,0},z_{2,0},...,z_{i,0}),\quad
\xi_{i}=(z_{i},...,z_{n}),\end{equation}where
$z_{1,0},z_{2,0},...,z_{n,0}$ are fixed points. Using Theorem 2.1,
Proposition 2.1 and Proposition 3.2 we obtain the main
result in this section.\\
\textbf{Theorem 3.1.} \emph{The $n{\times}n$ matrix function}
\begin{equation}W(\emph{z})=W_{1}(z)W_{2}(z_{1},\xi_{2}){\times}...{\times}W_{n}(\mu_{n-1},z_{n})
\end{equation}\emph{is a fundamental rational solution of KZ system (1.1). Here}
\begin{equation}\frac{\partial}{{\partial}z_{j}}W_{i}(\mu_{i-1},\xi_{i})=
{\rho}A_{i}(\mu_{i-1},\xi_{i})W_{i}(\mu_{i-1},\xi_{i})
\end{equation}
\emph{and}
\begin{equation}W_{i}(\mu_{i-1},\xi_{i})=I_{n}.\end{equation}
\section{KZ equations in terms of new variables}
Following A.Varchenko [13] we change the variables
\begin{equation}u_{1}=z_{1}-z_{2},\quad
u_{k}=\frac{z_{k}-z_{k+1}}{z_{k-1}-z_{k}},\quad
(2{\leq}k{\leq}n-1),\end{equation} \begin{equation}
u_{n}=z_{1}+z_{2}+...+z_{n}.\end{equation} The KZ system takes the
following form
\begin{equation}\frac{{\partial}W}{{\partial}u_{j}}=
{\rho}H_{j}(u)W,\quad (1{\leq}j{\leq}n),\end{equation}where
$u=(u_{1},u_{2},...,u_{n}).$ In this section we are going to
investigate the form of $H_{j}(u)$ in a more detailed manner than it
was done in paper [13]. We represent relations (4.1) and (4.2) in
the form
\begin{equation}z_{k}-z_{k+1}=u_{1}{\cdot}u_{2}{\cdot}...{\cdot}u_{k},\quad
(1{\leq}k{\leq}n-1),\quad  z_{1}+z_{2}+...+z_{n}=u_{n}.
\end{equation}  We write equality (4.4)in the matrix form
\begin{equation}SZ=U,\end{equation} where
\begin{equation} Z=\mathrm{col}[z_{1},z_{2},...,z_{n}],\quad
U=\mathrm{col}[u_{1},u_{1}{\cdot}u_{2},...,u_{1}{\cdot}u_{2}{\cdot}...{\cdot}u_{n-1},u_{n}].\end{equation}
The elements $s_{k,\ell}$ of the $n{\times}n$ matrix $S$ are equal
to zero except
\begin{equation}s_{k,k}=1,\quad s_{k,k+1}=-1,\quad
s_{n,k}=1.\end{equation} It is easy to check that the matrix
$S^{-1}$ has the following form
\begin{equation}S^{-1}=\frac{1}{n}\left[\begin{array}{cccccc}
                                          n-1 & n-2 & n-3 & ... & 1 & 1 \\
                                          -1 & n-2& n-3 & ...& 1 & 1\\
                                          -1 & -2 & n-3 & ... & 1 & 1 \\
                                          -1 & -2& -3& ...& 1 & 1 \\
                                          ...& ... & ...& ... & ... & ... \\
                                          -1& -2& -3 & ...& 1 & 1\\
                                          -1& -2 & -3 & ...& -n+1& 1
                                        \end{array}\right].\end{equation}
Using the formula
\begin{equation}\frac{{\partial}W}{{\partial}u_{j}}=
\sum_{k=1}^{n}\frac{{\partial}z_{k}}{{\partial}u_{j}}\frac{{\partial}W}{{\partial}z_{k}}
\end{equation}we have
\begin{equation}H_{j}(u)=\sum_{k=1}^{n}\frac{{\partial}z_{k}}{{\partial}u_{j}}A_{k}(u).\end{equation}
We note that ${\partial}z_{k}/{\partial}u_{j}$ are defined by the
relation
\begin{equation}\frac{{\partial}Z}{{\partial}u_{j}}=S^{-1}\frac{{\partial}U}{{\partial}u_{j}}.
\end{equation}It follows from (4.8) and (4.11) that
\begin{equation}\frac{{\partial}Z}{{\partial}u_{n}}=\frac{1}{n}\mathrm{col}[1,1,...,1].\end{equation}
In view of (1.2) and (4.12) the relation
\begin{equation}H_{n}(u)=\sum_{k=1}^{n}A_{k}(u)=0 \end{equation}
is true. It means that
\begin{equation}\frac{{\partial}W}{{\partial}u_{n}}=0.\end{equation}\emph{The
last relation is well-known} (see [3],[13]).\\
We introduce the triangular $n{\times}n$ matrix
\begin{equation}T=\left[\begin{array}{cccc}
                          1 & 1 & ... & 1 \\
                          0& 1 & ... & 1 \\
                          ... & ... & ... & ... \\
                          0 & 0 & ... & 1
                        \end{array}\right]\end{equation}
and the vector
\begin{equation}Y=TU=\mathrm{col}[y_{1},y_{2},...,y_{n}].\end{equation}
The matrices $T$ and $S$ are connected by the equality
\begin{equation}S^{-1}=T-\frac{1}{n}C,\end{equation} where
\begin{equation}C=\left[\begin{array}{ccccc}
                          1 & 2& ... & n-1 & n-1 \\
                          1& 2 & ... & n-1 & n-1\\
                          ... & ... & ... & ...& ... \\
                          1& 2 & ... & n-1 & n-1
                        \end{array}\right]\end{equation}
In view of (4.17) and (4.18) the equality (4.10) can be written in
the form
\begin{equation}H_{j}(u)=\sum_{k=1}^{n}\frac{{\partial}y_{k}}{{\partial}u_{j}}A_{k}(u).\end{equation}
Formulas (4.4), (4.15), (4.16) and (4.19) imply the following assertion.\\
\textbf{Proposition 4.1}\emph{The matrix functions $H_{j}$ do not
depend on $u_{1}$, when $2{\leq}j{\leq}n.$}\\
Let us consider separately $H_{1}$. In this case we have
\begin{equation}\frac{{\partial}U}{{\partial}u_{1}}=
\mathrm{col}[v_{1},v_{2},...,v_{n}],\end{equation}where
\begin{equation}v_{1}=1; \quad v_{k}=u_{2}{\cdot}...{\cdot}u_{k},\quad
(2{\leq}k{\leq}n-1); \quad v_{n}=0.\end{equation}According to
(4.15), (4.16) and (4.20) the equalities
\begin{equation}\frac{{\partial}y_{k}}{{\partial}u_{1}}=\sum_{j=k}^{n-1}v_{j}
\quad (1{\leq}k{\leq}n-1);\quad
\frac{{\partial}y_{n}}{{\partial}u_{1}}=0\end{equation} are true.
Taking into account the relations
\begin{equation}\frac{{\partial}y_{k}}{{\partial}u_{1}}-\frac{{\partial}y_{j}}{{\partial}u_{1}}
=\frac{z_{k}-z_{j}}{u_{1}}, \quad k>j,\end{equation} we deduce that
\begin{equation}H_{1}=\sum_{k>j}P_{kj}/u_{1}.\end{equation}
Thus we have proved the assertion.\\
\textbf{Proposition 4.2} \emph{The matrix function  $\mathrm{H}_{1}$
in KZ system (1.1) depends only on $\mathrm{u}_{1}$ and has form}
(4.24).\\
Let us consider the case when $2{\leq}k{\leq}n.$ In this case we
have
\begin{equation}\frac{{\partial}U}{{\partial}u_{k}}=\mathrm{col}[v_{1,k},v_{2,k},...,v_{n,k}],\end{equation}
where
\begin{equation}v_{p,k}=0,\quad (1{\leq}p{\leq}k-1);\quad
v_{p,k}=u_{1}v_{p}/u_{k},\\ (k{\leq}p{\leq}n-1);\end{equation}
\begin{equation}v_{n,k}=0,\quad k<n.\end{equation}Using again
(4.15), (4.16) and (4.20)we receive the relations
\begin{equation}\frac{{\partial}y_{s}}{{\partial}u_{k}}=\sum_{p=max(k,s)}^{n-1}v_{p,k}
\quad (1{\leq}k{\leq}n-1);\quad
\frac{{\partial}y_{n}}{{\partial}u_{k}}=0,\quad
(1{\leq}k<n).\end{equation}We introduce the denotation
\begin{equation}\alpha_{s,j,k}=(\frac{{\partial}y_{s}}{{\partial}u_{k}}-\frac{{\partial}y_{j}}{{\partial}u_{k}})/
(z_{s}-z_{j})\quad s>j,\end{equation}Relations (4.4),(4.29), and
(4.28) (4.29) imply that
\begin{equation}\alpha_{s,j,k}=\sum_{p=max(j,k)}^{s-1}v_{p}/(u_{k}\sum_{p=j}^{s-1}v_{p}).
\end{equation}Hence the following formula
\begin{equation}H_{k}=\sum_{s>j}\alpha_{s,j,k}(u)P_{s,j}\end{equation}
is valid. From formulas (4.26) and (4.30) we infer that
\begin{equation}\alpha_{s,j,k}=1/u_{k},\quad
s>j{\geq}k,\end{equation}
\begin{equation}\alpha_{s,j,k}=1+o(1),\quad k=j+1,\quad
s{\geq}j+2,\end{equation}
\begin{equation}\alpha_{s,j,k}=o(1),\quad k>j+1.\end{equation}
The equality $g(u)=o(u)$ means that the function $g(u)$ is regular
in the neighborhood of $u=0$ and $g(0)=0$. Using formulas (4.31)
-(4.34) we
deduce the theorem.\\
\textbf{Theorem 4.1} \emph{The following relations are true}
\begin{equation}H_{1}={\Omega}_{1}/u_{1},\quad H_{n}(u)=0,\end{equation}
\begin{equation}H_{s}={\Omega}_{s}/u_{s}+P_{s-1}+o(u),\quad
2{\leq}s{\leq}n-1,
\end{equation}
\emph{where}
\begin{equation}P_{r}=\sum_{j>r}P_{j,r},\quad \Omega_{s}=P_{s}+P_{s+1}+...+P_{n}.\end{equation}
\textbf{Remark 4.1.}The formula $H_{n}(u)=0$  and the first term of
asymptotic (4.36) were found by A Varchenko [13].The second term of
asymptotic (4.36) and explicit formulas (4.24) and  (4.31) are new.\\
\textbf{Remark 4.2.} In this section we consider arbitrary $\rho$
(not only integer).\\
\textbf{Corollary 4.1.} \emph{The fundamental solution $W(u)$ of KZ
system (1.1) can be represented in the form}
\begin{equation}
W(u)=W_{1}(u_{1})W(u_{2},...,u_{n}),\end{equation}\emph{where}
$W_{1}(u_{1})$ and $W(u_{2},...,u_{n})$ are fundamental solutions of
the system \begin{equation}
\frac{dW_{1}}{du_{1}}={\rho}H_{1}(u_{1})W_{1}(u_{1})
\end{equation}
\emph{and the system}
\begin{equation}
\frac{{\partial}W_{j}}{{\partial}u_{j}}={\rho}H_{j}(u_{2},...,u_{n})W_{2}(u_{2},...,u_{n}),\quad
2{\leq}j{\leq}n-1.
\end{equation}\emph{respectively.}
\section{Eigenvalues and eigenvectors of $\Omega_{s}$}
According to (4.35) the matrices $\Omega_{s}$ are coefficients of
the main parts of the asymptotic $H_{s}$. In the next section we
show that the eigenvalues and eigenvectors of $\Omega_{s}$ define
the asymptotic of the solution $W(u)$ of the KZ system. Therefore in
this section we find the eigenvalues and eigenvectors of
$\Omega_{s}$ in an explicit form. To do it we ise the following
result.\\
\textbf{Theorem 5.1.}\emph{The matrices $\Omega_{s}$  have the
structure}
\begin{equation}\Omega_{s}=\left[\begin{array}{cc}
                                  N_{s-1}I_{s-1} & 0 \\
                                  0 & \omega_{s}
                                \end{array}\right],
                                (1{\leq}s{\leq}n-1),
\end{equation}\emph{where the $(n-s+1){\times}(n-s+1)$ matrix
$\omega_{s}$ and the number $N_{s}$ are defined by the relations}
\begin{equation}\omega_{s}=\left[\begin{array}{ccccc}
                                   N_{s} & 1 & 1 & ... & 1 \\
                                   1 & N_{s}  & 1 & ... & 1 \\
                                   1 & 1 & N_{s} & ... & 1 \\
                                   1 & 1 & 1 & ... & N_{s}
                                 \end{array}\right],\end{equation}
\begin{equation}N_{s}=(n-s)(n-s-1)/2,\quad 1{\leq}s{\leq}n-1.\end{equation}
\emph{Proof.} It follows from (4.37) that non-diagonal elements of
$\Omega_{s}$ and the right part of (5.1) coincide. According to
(5.1) and (5.2) we have
\begin{equation}P_{ij}\Omega_{s}P_{ij}=\Omega_{s}\quad if\quad
either \quad 1{\leq}i,j{\leq}s-1 \quad or \quad
s{\leq}i,j{\leq}n,\end{equation}The relations (4.37 )imply that
$\Omega_{s}$ consists of $(n-s)(n-s+1)/2$ addend of $P_{ij}$. It
means that the left upper blocks of $\Omega_{s}$ and the right side
of (5.1) have the same diagonal elements $N_{s-1}$. Using again
(4.37) we prove that all other diagonal elements of   $\Omega_{s}$
are
equal to $N_{s}$. The theorem is proved.\\
From the block structure of $\Omega_{s}$ we deduce the assertion.\\
\textbf{Corollary 5.1.} \emph{The matrices $\Omega_{s},\quad
(1{\leq}s{\leq}n-1)$ have the common system of the eigenvectors}
\begin{equation}v_{k}=\mathrm{col}[\overbrace{0,0,...,0,}^{n-k-1},-k,\overbrace{1,1,...,1}^{k}],\quad
(1{\leq}k{\leq}n-1),\end{equation}
\begin{equation}v_{n}=\mathrm{col}[\overbrace{1,1,...,1}^{n}].\end{equation}
\emph{The eigenvalues of $\Omega_{s}$ corresponding to the vectors
$v_{k}$ are defined by the relations}
\begin{equation}\lambda_{k,s}=N_{s-1},\quad
n{\geq}k>n-s+1,\end{equation}
\begin{equation}\lambda_{k,s}=N_{s}-1,\quad
1{\leq}k{\leq}n-s+1.\end{equation} \textbf{Corollary 5.2.} \emph{All
the eigenvalues of $\Omega_{s}$ are integer and non-negative.}
\section{Asymptotic of solutions of the KZ system}
By $D$ we denote the domain
\begin{equation}z_{1}>z_{2}>...>z_{n}.\end{equation} It means that
(see (4.1))
\begin{equation}u_{1}>u_{2}>...>u_{n}.\end{equation}The following
proposition was obtained by A. Varchenko [13].\\
\textbf{Theorem 6.1.} \emph{Let us assume that $\rho$ is irrational.
Then:\\
1) For every vector $v_{k}$ there exists a unique solution
$\psi_{k}$ of the KZ system in $D$ such that }
\begin{equation}\psi_{k}(u)=[\prod_{s=1}^{n-1}u_{s}^{\rho\lambda_{k,s}}](v_{k}+o(u)),\quad
(1{\leq}k{\leq}n),\end{equation} 2)\emph{The solutions $\psi_{k}(u)$
form a basis in the space of the solutions of the KZ system on $D$.}\\
We deduce the following addition to Varchenko Theorem.\\
\textbf{Theorem 6.2.} \emph{Let us consider  KZ system (1.1), where
$\rho$ is integer. The assertions 1) and 2) of Varchenko theorem 6.1
are true.}\\
\emph{Proof.} It follows from Theorem 3.1  that there exists the
fundamental rational solution of KZ system (1.1). Then according to
the results from ([8],section 2) this solution can be represented in
form (6.3). The
theorem is proved.\\
\textbf{Remark 6.1.} When $\rho$ is integer the representation (6.3)
is true not only in the domain $D$ but in a neighborhood of $u=0.$
\section{Examples}
\textbf{Example 7.1.} Let us consider the simplest case n=3. Using
(4.31), (4,35)  and (4.39), (4.40) we have
\begin{equation}\frac{{\partial}W}{{\partial}u_1}=\rho{\left(\frac{\Omega_1}{u_1}\right)}W,\end{equation}
\begin{equation}\frac{{\partial}W}{{\partial}u_2}=\rho{\left(\frac{P_{3,2}}{u_2}+
\frac{P_{3,1}}{1+u_2}\right)}W,\end{equation}where
\begin{equation}
\Omega_1=P_{1,2}+P_{1,3}+P_{23}=\left[\begin{array}{ccc}
                                                      1 & 1 & 1 \\
                                                      1 & 1& 1 \\
                                                      1 & 1 & 1
                                                   \end{array}\right].\end{equation}
According to Corollary 4.1 we can represent $W(u)$ in the form
\begin{equation}W(u)=W_{1}(u_{1})W_{2}(u_{2}),\end{equation} where
$W_{1}(u_{1})$ and $W_{2}(u_{2})$ are fundamental solutions of
system (7.1) and (7.2) respectively. It is easy to see that
\begin{equation}W_{1}(u_{1})=u_{1}^{\rho\Omega_1}.\end{equation}
We note that the eigenvalues of $\Omega_1$ are integers
($\lambda_1=3,\quad \lambda_2=\lambda_3=0$). Hence the matrix
$W_{1}(u_{1})$ is rational when $\rho$ is integer. We introduce the
matrix function
\begin{equation}F(y)=W_{2}(y)y^{\rho}(1+y)^{\rho},\quad
y=u_{2}.\end{equation}Relations (7.2) and (7.6) imply that
\begin{equation}\frac{dF}{dy}=\rho{\left(\frac{P_{3,2}+I}{y}+
\frac{P_{3,1}+I}{1+y}\right)}F.\end{equation}We consider the
constant vectors
\begin{equation}w_{1}=\mathrm{col}[0,1,-1],\quad
w_{2}=\mathrm{col}[1,-2,1].\end{equation}It is easy to see that
\begin{equation}(P_{3,2}+I)w_{1}=0,\quad
(P_{3,1}+I)w_{1}=-w_{2},\end{equation}
\begin{equation}(P_{3,2}+I)w_{2}=3w_{1}+2w_{2},\quad
(P_{3,1}+I)w_{2}=2w_{2},\end{equation} We represent $F(y)$ in the
form
\begin{equation}F(y)=\phi_{1}(y)w_{1}+\phi_{2}(y)w_{2}\end{equation}
and substitute it in (7.7). In view of (7.7) and (7.11) we have
\begin{equation}\phi_{1}^{\prime}=\frac{3\rho}{y}\phi_{2},\quad
\phi_{2}^{\prime}=\rho(\frac{2\phi_{2}}{y}+\frac{2\phi_{2}}{1+y}-\frac{\phi_{1}}{1+y}).
\end{equation}It follows from (7.12) that
\begin{equation}y(1+y)\phi_{1}^{\prime\prime}(y)+[1+y-2\rho(1+2y)]\phi_{1}^{\prime}
+3{\rho}^{2}\phi_{1}=0.\end{equation}By introducing
$\psi(y)=\phi_{1}(-y)$ we reduce equation (7.13) to Gauss
hypergeometric equation [1]:
\begin{equation}y(1-y)\psi^{\prime\prime}(y)+[\gamma-(\alpha+\beta+1)y]\psi^{\prime}(y)
-\alpha\beta\psi(y)=0,\end{equation}where
\begin{equation}\alpha=-\rho,\quad \beta=-3\rho,\quad
\gamma=1-2\rho.\end{equation}Let $\psi_{1}$ and $\psi_{2}$ be
linearly independent solutions of equation (7.14). Then the vector
functions
\begin{equation}Y_{k}(u_{2})=[\psi_{k}(-u_{2})w_{1}-\frac{3\rho}{u_{2}}\psi_{k}^{\prime}(-u_{2})w_{2}]
u_{2}^{-\rho}(1+u_{2})^{-\rho},\quad k=1,2\end{equation} are the
solutions of system (7.2). It is easy to check that the vector
function
\begin{equation}Y_{3}(u_{2})=u_{2}^{\rho}(1+u_{2})^{\rho}\mathrm{col}[1,1,1]\end{equation}
is the solution of system (7.2) too. Hence we deduced the assertion\\
\textbf{Proposition 7.1.} \emph{The fundamental solution $W(u)$ of
system (7.1), (7.2) is defined by relations (7.4) and (7.5) , where}
\begin{equation}W_{2}(u_{2})=[Y_{1}(u_{2}),Y_{2}(u_{2}),Y_{3}(u_{2})].\end{equation}
\textbf{Remark 7.1.} Substitution (7.11) was prompted by book ([3],
section 4.2).\\
As a by-product we deduce from Proposition 7.1 the following
result.\\
\textbf{Proposition 7.2.} \emph{The solutions of Gauss
hypergeometric equation (7.14) are rational functions if
\begin{equation}\alpha=-\rho,\quad \beta=-3\rho,\quad
\gamma=1-2\rho,\end{equation}where $\rho$ is integer.}\\
\textbf{Remark 7.2.}The KZ system is connected with the generalized
hypergeometric equations of several variables (see [5]). Hence the
Theorem 2.1 can be applied to the generalized hypergeometric
equations of several variables .
 \textbf{Example 7.2.} Let us consider
separately the case $n=3,\rho=-1.$ Using results of paper [7] we
have
\begin{equation}\phi_{1}=\frac{1}{1-y},\quad
\phi_{2}=\frac{1}{y^{2}}-\frac{1}{y}.\end{equation}The corresponding
hypergeometric equation (7.14) has the form\
\begin{equation}y(1-y)\psi^{\prime\prime}(y)+(3-5y)\psi^{\prime}(y)
-3\psi(y)=0,\end{equation} \textbf{Example 7.3.} Now we consider the
case n=4. We have
\begin{equation}\frac{{\partial}W}{{\partial}u_1}=\rho{\left(\frac{\Omega_1}{u_1}\right)}W,\end{equation}
\begin{equation}\frac{{\partial}W}{{\partial}u_2}=\rho{\left(\frac{\Omega_{2}}{u_2}+
\frac{P_{1,3}}{1+u_2}+\frac{P_{1,4}(1+u_{3})}{1+u_2+u_{2}u_{3}}\right)}W,\end{equation}
\begin{equation}\frac{{\partial}W}{{\partial}u_3}=\rho{\left(\frac{P_{4,3}}{u_3}+
\frac{P_{4,2}}{1+u_3}+\frac{P_{4,1}u_{2}}{1+u_2+u_{2}u_{3}}\right)}W,\end{equation}where
\begin{equation}\Omega_{1}=\left[\begin{array}{cccc}
                                   3 & 1 & 1 & 1 \\
                                   1 & 3 & 1 & 1 \\
                                   1 & 1 & 3 & 1 \\
                                   1 & 1 & 1 & 3
                                 \end{array}\right],\end{equation}
                                 \begin{equation}
                                   \Omega_{2}=\left[\begin{array}{cccc}
                                   3 & 0 & 0 & 0 \\
                                   0 & 1 & 1 & 1 \\
                                   0& 1 & 1& 1 \\
                                   0 & 1 & 1 & 1
                                 \end{array}\right].\end{equation}

\textbf{Remark 7.3.} For the case $n=4,\quad \rho=-1$ we constructed
the solution $W(z)$ of system (2.4) in the explicit form. Using
formula (3.17) we can obtain in the explicit form the solution of KZ
system
(1.1) when $n=4,\quad \rho=-1.$ \\
\textbf{Acknowledgement.}  I am very grateful to A. Tydnyuk. His
calculations [11],[12]  were very helpful for my work.
\begin{center}\textbf{References} \end{center}
1. Bateman H. and Erdely A.,Higher Transcendental Functions, v.1,New
York, 1953\\
2. Burrow M., Representation Theory of Finite Groups, Academic
Press,
1965.\\
3. Etingof P.I., Frenkel I.B., Kirillov A.A. (jr.), Lectures on
Representation Theory and Knizhnik-Zamolodchikov Equations, Amer.
Math. Society, 1998.\\
4. Felder G. and Veselov A., Polynomial Solutions of the
Knizhnik-Zamolod-\\chikov Equations and Schur-Weyl Duality,
International Mathematics Research Notices, v.2007,
 2007,1-27.\\
5.  Matsuo A., An Applications of Aumoto-Gelfand Hypergeometric
Functions to the SU(n) Knizhnik-Zamolodchikov Equations, Comm.in
Math.Phys. v.134, Number 1,65-77, 1990.\\
6. Sakhnovich L.A., Spectral Theory of Canonical Differential
Systems.
Method of Operator Identities., OT,Adv. and Appl.,v.107, 1999.\\
7. Sakhnovich L.A., Rational Solutions of KZ Equation, Existence and
Construction, arXiv.math:-ph/0609067, 2006.\\
8. Sakhnovich L.A., Meromorphic Solutions of Linear Differential
Systems, Painleve Type Functions,Operator and Matrices, v.1,\\
Number 1, 87-111, 2007.\\
 9. Sakhnovich L.A., Rational Solutions of
KZ Equation, Case $S_{4}$,\\ arXiv:math.CA/0702404, 2007.\\ 10.
Sakhnovich L.A.,Explicit  Rational Solutions of
 Knizhnik-Zamolodchikov
Equation, Central European Journal of Math., v.1, Number 1,
179-187,2008.\\ 11. Tydnyuk A., Rational Solution of KZ Equation
(Example),\\ arXiv:math/0162153,1-13, 2006.\\ 12. Tydnyuk A.,
Explicit Rational Solution of KZ Equation,\\ arXiv:math/07091141,
2007.\\ 13. Varchenko A., Asymptotic Solutions of
Knizhnik-Zamolodchikov Equation and Crystal Base, Comm.Math.Phys.
171, 99-138,1995.\\
\end{document}